\documentclass{amsproc}

\newtheorem{Theorem}{Theorem}

\newtheorem{theorem}{Theorem}[section]
\newtheorem{lemma}[theorem]{Lemma}
\newtheorem{proposition}[theorem]{Proposition}

\theoremstyle{definition}

\theoremstyle{remark}

\numberwithin{equation}{section}

\newcommand{\abs}[1]{\lvert#1\rvert}

\newcommand{\R}{{\mathbb R}}
\newcommand{\Z}{{\mathbb Z}}

\newcommand{\A}{{\mathbb A}}
\newcommand{\N}{{\mathbb N}}
\newcommand{\LL}{{\mathbb L}}
\newcommand{\IA}{{\stackrel{\circ}{\A}}}
\newcommand{\PP}{{\mathcal P}_f(\A)}
\begin{document}

\title[NICIS homeomorphisms]{Annular area preserving
homeomorphisms which admit no interior compact invariant sets}

\author{Shigenori Matsumoto}
\address{Department of Mathematics, College of
Science and Technology, Nihon University, 1-8-14 Kanda, Surugadai,
Chiyoda-ku, Tokyo, 101-8308 Japan
}
\email{matsumo@math.cst.nihon-u.ac.jp
}
\thanks{The author is partially supported by Grant-in-Aid for
Scientific Research (C) No.\ 20540096.}


\date{\today }

\begin{abstract}
For any irrational number $\alpha$, there exists an ergodic area
 preserving homeomorphism of the closed annulus which is isotopic
to the identitity, admits no compact invariant set contained in the
interior of the annulus, and has the rotation number $\alpha$.
\end{abstract}

\maketitle

\section{Introduction}
Let $\A=S^1\times[-1,1]$ be the closed annulus equipped with
the standard area. Denote by $G^r$ ($r=0,1,\cdots,\infty$)
the group of the area preserving
$C^r$ diffeomorphisms of $\A$ which are isotopic to the identity,
and by $\partial_\pm\A=S^1\times\{\pm 1\}$ the upper or lower boundary
curves of the annulus.
In \cite{He1}, Michael Hermann announced that if $f\in G^\infty$
has no periodic points and if the rotation number of $f$ restricted
to $\partial_-\A$ is a non Liouville number, then $f$ is $C^\infty$
conjugate to a rigid rotation near $\partial_-\A$.
According to \cite{FS}, in fact he seems to have shown the following.

{\em If the rotation number of $f\vert_{\partial_-\A}$ is non Liouville,
then there is a family of invariant smooth circles accumulating to
$\partial_-\A$.}

B. Fayad and M. Saprykina \cite{FS} then showed the optimality of this result
by proving that for any Liouville number
$\alpha$, there is $f\in G^\infty$ such that $f$ is weakly mixing
and that the rotation number of the two boundary curves is $\alpha$.

In this paper we consider a somewhat different 
topological counter part of the Herman
phenomenon. A map in $G^r$ is called a {\em NICIS map} if it admits no
compact invariant set contained in the interior of $\A$.

Let $\tilde\A=\R\times[-1,1]$ be the universal cover of $\A$ and
$p:\tilde A\to\R$ be the projection onto the first factor.
Fix once and for all a lift $\tilde f:\tilde\A\to\tilde\A$ of $f\in
G^r$.
Then the map $p\circ\tilde f-p:\tilde\A\to\R$ is invariant by the covering
transformation group and defines a map on $\A$, which we shall denote by
the same letter $p\circ\tilde f-p$. For a NICIS map $f\in G^0$,
there is an irrational number $\alpha$ such that for any
$f$-invariant probability measure $\mu$ on $\A$, we have
$$
\langle p\circ\tilde f-p,\mu\rangle=\alpha.$$
See Appendix A. The number $\alpha$ is called the {\em rotation number}
of $\tilde f$ and is denoted by $\rho(\tilde f)$. 
The main result of the present paper is the following.

\begin{Theorem} \label{T}
(1) If $\alpha$ is a Liouville number, then there is a NICIS $C^\infty$-
diffeomorphism $f$ (and its lift $\tilde f$) such that $\rho(\tilde f)=\alpha$.

(2) For any irrational number $\alpha$, 
there is a NICIS homeomorphism $f$ such that $\rho(\tilde f)=\alpha$. 
\end{Theorem}

Recall that an irrational number $\alpha$ is called a {\em Liouville
number}
if for any $\tau>0$ and $\varepsilon>0$, there is a rational number 
$p/q$, $(p,q)=1$, such that $\abs{\alpha-p/q}<\varepsilon/q^\tau$.

The proof of (1), given in Sect.\ 2, is based on
 the fast approximation by conjugacy method developed in \cite{FS}.
But since the NICIS diffeomorphism is much simpler to handle than the
weakly mixing diffeomorphism treated in \cite{FS}, 
our argument becomes vastly easier than in \cite{FS}.
However, as will be apparent in Sect.\ 2,
it looks hopeless to prove (2) by this method,
and we show it in Sect.\ 3 by constructing a conjagacy with a Anzai skew product 
$F:S^1\times\R\to S^1\times\R$ over the rotation $R_\alpha:S^1\to S^1$
by $\alpha$.  The obtained homeomorphism is not weakly mixing,
although it is ergodic.
In Appendix A we show some fundamental properties of
NICIS maps and in Appendix B we prepare necessary prerequisites about
Anzai skew products.

\section{Fast approximation by conjugation method}

Define the $C^r$-distance $d_r$ on $G^r$ by
$$
d_r(f,g)=\max\{\Vert f-g\Vert_r,\Vert f^{-1}-g^{-1}\Vert_r\},
$$
where $\Vert\cdot\Vert_r$ denotes the $C^r$-norm.
Let us denote the standard $S^1$-action on $\A$ by $S_t$ ($t\in S^1$):
$$
S_t(x,y)=(x+t,y).$$

For any $n\in\N$, choose an increasing sequence of
 positive numbers $a_n$ tending to $1$.
Let $$U_n=S^1\times([-1,-a_n)\cup(a_n,1])$$ and
$$G^r_n=\{f\in G^r\mid O_f(x)\cap U_n\neq \emptyset,\  \forall x\in \A\},$$
where
$O_f(x)$ denotes the $f$-orbit of $x$.
Notice that $G^r_n$ is $C^0$-open in $G^r$
and that $f\in G^r$ is NICIS if and only if $f\in\cap_n G^r_n$.

For the moment let $\alpha$ be an arbitrary irrational number.
We shall construct a diffeomorphism $h_n\in G^\infty$ and
a rational number $\alpha_n=p_n/q_n$ ($n\in\N$)
such that 

\smallskip
\noindent
(A) $h_n\circ S_{\alpha_n}=S_{\alpha_n}\circ h_n$.

\smallskip
Denoting
$$
H_n=h_1\circ h_2\circ\cdots\circ h_n,\ \ 
f_n=H_n\circ S_{\alpha_{n+1}}\circ H_n^{-1},
$$
we will also require the following.

\smallskip
\noindent
(B) {\em The map $f_n$ converges to some $f\in G^\infty$ in the
$C^\infty$-topology
and
$\alpha_n\to\alpha$.}

\smallskip
\noindent
(C) {\em For any $n$, the $C^0$-closure of the set
$\{f_m\vert m\geq n\}$ is contained in $G^0_n$.}

\smallskip
Then the limit $f$ is a NICIS $C^\infty$-diffeomorphism and by the
continuity of the rotation number we have 
$\rho(\tilde f)=\alpha$ for some lift $\tilde f$.
The condition (A) plays a key role in getting (B) and (C), since it
implies the
equalities
$$
f_{n+1}-f_n=H_{n+1}\circ S_{\alpha_{n+2}}\circ H_{n+1}^{-1}
-H_{n+1}\circ S_{\alpha_{n+1}}\circ H_{n+1}^{-1},
$$
$$
f_{n+1}^{-1}-f_n^{-1}=H_{n+1}\circ S_{-\alpha_{n+2}}\circ
H_{n+1}^{-1}
-H_{n+1}\circ S_{-\alpha_{n+1}}\circ
H_{n+1}^{-1},$$
which enable us to estimate $d_r(f_{n+1},f_n)$
for example by Lemma \ref{L} below.

Assume we already defined $h_1,\cdots,h_{n-1}$ and
$\alpha_1,\cdots,\alpha_n=p_n/q_n$.
The construction of $h_n$ is carried out in the following way.
We identify the cyclic $q_n$-covering space of $\A$
with $\A$ in a standard way.
First we construct $k_n\in G^\infty$ which keeps the
vertical line $I_0=0\times[-1,1]$ invariant, and let 
$h_n$  be the lift of $k_n$ which keeps $I_0$ invariant.
Then clearly $h_n$ commutes with $S_{1/q_n}$,
and hence with $S_{\alpha_n}$.

The construction of $k_n$ goes as follows. Fix real numbers
$b_n$ and $c_n$ so that $a_n<b_n<c_n<a_{n+1}$.
Choose $k_n\in G^\infty$ such that $k_n$ is the identity on
the neighbourhood $U_{n+1}$ and on the vertical line $I_0$
 and that $k_n$ maps the horizontal curve
$S^1\times 0$ to the graph of the function $y=c_n\sin2\pi x$.
The Moser lemma
(\cite{M}) shows the existence of such an area preserving diffeomorphism.

Notice that the orbit foliation of the $S^1$-action $S_t$
is mapped by $k_n$ to a foliation all of whose leaves intersect
$S^1\times([-1,-c_n]\cup[c_n,1])$. 
Clearly this property is inherited to $h_n$ and also to
$H_n=h_1\circ \cdots h_{n-1}\circ h_n$ 
since $h_1\circ\cdots h_{n-1}$ is the identity
on $U_n$. That is,
any orbit of the $S^1$-action $H_n\circ S_t\circ H_n^{-1}$ ($t\in S^1$)
intersects $S^1\times([-1,-c_n]\cup[c_n,1])$.

Next choose a rational number $\alpha_{n+1}=p_{n+1}/q_{n+1}$ which approximates
the given $\alpha$ so as to satisfy the following conditions (2.1), (2.2)
and (2.3) below.

\medskip
\noindent
(2.1) {\em The number $q_{n+1}$ is big enough so that any orbit of the 
diffeomorphism
$H_n\circ S_{1/q_{n+1}}\circ H_n^{-1}$ intersects 
$S^1\times([-1,-b_n]\cup[b_n,1])$. We also assume
$$
\abs{\alpha-\alpha_{n+1}}<\abs{\alpha-\alpha_n}$$
for any $n$.}

\medskip
The same intersection property is enjoyed also by the 
diffeomorphism
$f_n=H_n\circ S_{\alpha_{n+1}}\circ H_n^{-1}$, 
since the both diffeomorphims, periodic of
order $q_{n+1}$, have the same orbit structure.

Furthermore if $g\in G^0$ satisfies
$$\Vert g^j-f_n^j\Vert_0<b_n-a_n\ \ \ (1\leq j\leq q_{n+1}),$$
then 
$g\in G^0_n$, i.\ e.\ any orbit of $g$ intersects $U_n$.
Let
$$
\varepsilon_n=\min\{\frac{b_n-a_n}{2},\frac{b_{n-1}-a_{n-1}}{2^2},\cdots
\frac{b_1-a_1}{2^{n}}\}.$$
If we choose diffeomorphisms $f_n=H_n\circ S_{\alpha_{n+1}}\circ H_n^{-1}$
in such a way that 
$$
\Vert f_{n+1}^j-f_n^j\Vert_0<\varepsilon_n\ \ (1\leq j\leq q_{n+1})
$$
for any $n$, then the sequence $\{f_n\}$ satisfies (C),
 and if further $f_n$ converges to some
$f$ in $G^0$, then $f$ belongs to $\cap_n G_n^0$,
that is, $f$ is a NICIS map.

Denoting the Lipschitz constant by $L(\cdot)$, we have
$$
\Vert f^j_{n+1}-f^j_n\Vert_0=\Vert H_{n+1}\circ S_{j\alpha_{n+2}}
\circ H_{n+1}^{-1}-H_{n+1}\circ S_{j\alpha_{n+1}}\circ H_{n+1}^{-1}\Vert_0\\
\leq L(H_{n+1})j\abs{\alpha_{n+2}-\alpha_{n+1}}.
$$
Since $L(h_n)\leq L(k_n)q_n$, $1\leq j\leq q_{n+1}$ and
$$
\abs{\alpha_{n+2}-\alpha_{n+1}}<2\abs{\alpha-\alpha_{n+1}}
$$
by (2.1), we get
$$
\Vert f^j_{n+1}-f^j_n\Vert_0\leq C_nq_1\cdots q_{n}q_{n+1}^2
2\abs{\alpha-\alpha_{n+1}},
$$
where $C_n=L(k_1)\cdots L(k_n)$.

Thus the following condition (2.2) assures 
$\Vert f_{n+1}^j-f_n^j\Vert_0<\varepsilon_n$.

\medskip
\noindent
(2.2) {\em We have}
$$
\abs{\alpha-\frac{p_{n+1}}{q_{n+1}}}<\frac{\varepsilon_n}{2C_nq_1
\cdots q_n}\cdot\frac{1} {q_{n+1}^2}.
$$

\medskip
Finally to guarantee the convergence of $f_n$, we need the following lemma.
The statement (1) is just Lemma 5.6 of \cite{FS}, and (2) is easy to
establish.

\begin{lemma} \label{L}
(1) For $r>0$, $\alpha,\beta\in\R$ and for $H\in G^r$,
we have
$$
d_r(H\circ S_\alpha\circ H^{-1}, H \circ S_\beta\circ H^{-1})
\leq C_2(r)\Vert H \Vert_{r+1}^{r+1}\abs{\alpha-\beta},
$$
where $C_2(r)$ is a constant depending only on $r$.

(2) We have
 $$\Vert h_n\Vert_r \leq C_3(n,r)q_n^{r},$$
 where 
$C_3(n,r)$ is a constant depending only on $n$ and $r$.
\qed
\end{lemma}

The lemma implies that
$$
d_r(f_{n+1},f_n)\leq C_4(n,r)q_1^{(r+1)^2}\cdots q_n^{(r+1)^2}q_{n+1}^{(r+1)^2}
\abs{\alpha-\alpha_{n+1}},
$$
where $C_4(n,r)$  is a constant depending only on $n$ and $r$.

Therefore the following (2.3) garantees that 
$d_n(f_{n+1},f_n)\leq 1/2^n$ and that $f_n$ converges in $G^\infty$.

\medskip
\noindent
(2.3) {\em  We have}
$$
\abs{\alpha-\frac{p_{n+1}}{q_{n+1}}}\leq \frac{1}
{2^nC_4(n,n)q_1^{(n+1)^2}\cdots q_n^{(n+1)^2}}\cdot\frac{1}{q_{n+1}^{(n+1)^2}}.
$$

\medskip
Now if $\alpha$ is a Liouville number, then it is possible to choose
$p_n/q_n$ successively so as to satisfy (2.1), (2.2) and (2.3) above.
This shows (1) of Theorem \ref{T}.

Notice that if we are working in $G^0$, the condition for 
$C^0$-convergence becomes;

\medskip
\noindent
(2.3')
$$
\abs{\alpha-\frac{p_{n+1}}{q_{n+1}}}\leq \frac{1}
{2^nC_4(n,0)q_1\cdots q_n}\cdot\frac{1}{q_{n+1}}.
$$
 and any irrational number $\alpha$ admits a sequence
$\{p_n/q_n\}$ satisfying (2.1) and (2.3').
On the other hand the condition (2.2) cannot be relaxed,
and it remains still necessary to require the well approximability of $\alpha$.
Recall an irrational number $\alpha$ is called 
{\em well approximable} if for any $\varepsilon>0$, there is a rational number 
$p/q$ such that $\abs{\alpha-p/q}<\varepsilon/q^2$.

\section{Anzai skew products}

The purpose of this section is to give a proof of Theorem \ref{T} (2).
Denote by $R_\alpha:S^1\to S^1$ the rotation by $\alpha$, and let $\varphi$ be a real valued continuous function on $S^1$.
Define a homeomorphism $F_{\alpha,\varphi}$ of the annulus $S^1\times\R$ by
$$
F_{\alpha,\varphi}(x,y)=(R_{\alpha}(x),y+\varphi(x)).
$$
The homeomorphism $F_{\alpha,\varphi}$ is called a {\em Anzai skew product}
over $R_\alpha$. Throughout this section we assume

\smallskip
\noindent
(3.1) {\em The rotation number $\alpha$ is irrational.}

\smallskip
\noindent
(3.2) {\em The function $\varphi$ satisfies}
$
\int_{S^1}\varphi(x)dx=0.
$

\smallskip
\noindent
(3.3) {\em The function
$\varphi$ is nonintegrable, i.\ e.\  the functional equation
$h\circ R_\alpha-h=\varphi$ has no continuous solution $h:S^1\to\R$.}

\smallskip
\noindent
(3.4) {\em The function $\varphi$ is absolutely continuous.}

\medskip

Notice that if condition (3.2) is not fulfilled, then 
the unique ergodicity of $R_\alpha$ implies that
all the orbits of $F_{\alpha,\varphi}$ are discrete in $S^1\times\R$
and there is no $F_{\alpha,\varphi}$-invariant probability measure
on $S^1\times\R$.
Also if there is a solution $h$ in (3.3), then the graph
of $h$ as well as its vertical translates will be an invariant curve. On the other hand
conditions (3.3) implies that there are no compact invariant set of
$F_{\alpha,\varphi}$ (Theorem 14.11, \cite{GH}).

For a point $(x,y)\in S^1\times\R$, denote by $\lambda(x,y)$ the limit set
of $(x,y)$ by the homeomorphism $F_{\alpha,\varphi}$,
 {\em i.\ e.} the union of the $\alpha$-limit set and the
$\omega$-limit set:
$$
\lambda(x,y)=\alpha(x,y)\cup\omega(x,y).
$$
Following H. Poincar\'e, define
$$
\LL_x=\{y\in\R\mid (x,y)\in\lambda(x,0)\}.
$$

Since $F_{\alpha,\varphi}$ commutes with the vertical translations,
$\LL_x$ forms a closed subsemigroup of $\R$, and
we have
$$
\lambda(x,y)\cap(x\times\R)=(x,y+\LL_x).
$$
Define
$$
D=\{x\in S^1\mid\LL_x=\R\},
$$
and 
notice that $x\in D$ if and only if the
orbit of $(x,y)$ is dense in $S^1\times\R$ for some (any) $y\in\R$.
Also define
$$
P=\{x\in S^1\mid\LL_x=[0,\infty)\},
\ \ 
N=\{x\in S^1\mid\LL_x=(-\infty,0]\}.
$$

The following theorem plays a crucial role in what follows.
It was proven in \cite{K} for $C^1$-function $\varphi$ 
in the situation where the $\lambda$-limit set
is replaced by  the 
$\omega$-limit set. But this is not enough
for our purpose. We shall include the proof in Appendix B.

\begin{theorem} \label{t4}
Assume that $F_{\alpha,\varphi}$ satisfies {\rm (3.1) $\sim$ (3.4)}.
Then all the three subsets $D$, $P$ and $N$ are nonempty and they
constitute a partition of $S^1$ into Borel sets.
\end{theorem}

Here is an outline of the proof of Theorem \ref{T} (2). Given any
irrational number $\alpha$, we shall construct a function
$\varphi:S^1\to\R$ satisfying (3.2) $\sim$ (3.4).
Moreover $\varphi$ is to satisfy (3.5) and (3.6) below.

\smallskip
\noindent
(3.5) {\em There exists a function $h\in L^1(S^1,dx)$ such that
$h\circ R_\alpha-h=\varphi$, where $dx$ denotes the Lebesgue measure on
$S^1$.}

\smallskip
\noindent
(3.6) {\em We have $\varphi(-x-\alpha)=\varphi(x)$ for any $x\in S^1$.}

\smallskip
We identify $S^1\times\R$ with the interior $\IA$ of the closed
annulus $\A$ in a standard way by
a diffeomorphism which maps a fiber $x\times\R$ onto a fiber
$x\times(-1,1)$. 
Clearly the map $F_{\alpha,
\varphi}$ extends to a homeomorphism of the closed annulus $\A$
in such a way that the restriction to each boundary component
is the rotation by $\alpha$.
Using a solution $h$ of (3.5), let us define a measurable map
$\Gamma(h):S^1\to\A$ by $\Gamma(h)(x)=(x,h(x))$.
Then we get a probability measure $\mu=\Gamma(h)_*dx$ on $\A$.
Condition (3.5) implies that
$\mu$ is invariant by $F_{\alpha,\varphi}$. 

Clearly $\mu$ is nonatomic, and $\mu(\partial \A)=0$. 
Furthermore once we show that $\mu$ is of full support,
we can find a homeomorphism $g$ of $\A$, identity on the boundary,
such that $g_*\mu$ is the standard area of $\A$ (\cite{OU}). 
Then the conjugate $g\circ F_{\alpha,\varphi}\circ g^{-1}$ is
a NICIS homeomorphism with rotation number equal to $\alpha$.
Also it is ergodic since it is metrically isomorphic to the rotation
$R_\alpha$ on $S^1$.

\smallskip
Given $\alpha$, the function $\varphi$ is constructed as follows.
For $x\in \R$, $\Vert x\Vert$ denotes the distance in $S^1=\R/\Z$ of 
the class of $x$ to $0$.
Choose a sequence $\{q_n\}_{n=1}^\infty$ of positive integers satisfying
$$
(3.7)  \ \Vert q_n\alpha\Vert <1/ q_n, \ \
(3.8) \ q_{n+1}> 10q_n,
$$
and define
$$
\varphi_N(x)=\sum_{n=1}^N
\frac{1}{ni}((e^{2\pi i q_n\alpha} -1)e^{2\pi i q_n x}
             -(e^{-2\pi i q_n\alpha} -1)e^{-2\pi i q_n x}).
$$
We have by (3.7)
$$
\vert e^{\pm 2\pi i q_n\alpha} -1\vert < 2\pi\Vert q_n\alpha\Vert<
\frac{2\pi}{q_n}.
$$
This, together with the assumption (3.8) shows that $\varphi_N$
converges  uniformly to a continuous function
 $\varphi$. To show that $\varphi$ is 
absolutely continuous, we have;
$$
\int_{S^1}\vert \varphi'(x)\vert dx
\leq(\int_{S^1}\vert \varphi'(x)\vert^2dx)^{1/2}
$$
$$
\leq(2\sum_{n=1}^\infty \frac{1}{n^2}\vert e^{2\pi i q_n\alpha}-1
\vert^2(2\pi q_n)^2)^{1/2}
\leq 4\sqrt{2}\pi^2\sum_n\frac{1}{n^2} <\infty.
$$
Since $\varphi'\in L^1(S^1,dx)$,
the function $\varphi$ is absolutely continuous.

The Fourier series of
the solution $h$ of the equation $h\circ R_\alpha-h=\varphi$ is given
by 
$$
h(x)=\sum_{n=1}^\infty \frac{1}{ni}(e^{2\pi i q_nx}-e^{-2\pi iq_nx})
=\sum_{n=1}^\infty \frac{1}{2n}\sin 2\pi q_n x.
$$
Direct computation shows that $h\in L^2(S^1,dx)\subset L^1(S^1,dx)$.
On the other hand $h$ is not continuous. 
To see this, recall that the continuity of $h$ is equivalent to
the uniform convergence of the associated Ces\'aro sum 
(Theorem 73 of \cite{HR}).  Let
$$C_n=\{x\in S^1\mid \sin2\pi q_nx\geq 1/2\}.
$$
Then by (3.8) a connected component of $C_n$ contains a 
connected component of $C_{n+1}$, and thus one can choose a point $x_0$
from the intersection of all the $C_n$'s. Then clearly the Ces\'aro
sum diverges at $x_0$.

The proof that the function $\varphi$ 
satisfies the conditions (3.2) $\sim$ (3.5)
is now complete, while the condition (3.6) can be verified
by a direct computation.

\medskip

Finally let us show that the measure $\mu$ has full support.
Notice that condition (3.6) implies that 
$$
J\circ F_{\alpha,\varphi}\circ J=F^{-1}_{\alpha,\varphi},$$
where $J$ is the involution of $S^1\times\R$ defined
by
$$
J(x,y)=(-x,-y)).
$$

Since our definition of $\LL_x$ is invariant by time reversion,
this shows that $-P=N$ on $S^1$. Therefore the measures of $P$
and $N$ must be the same. But since $R_\alpha$ is ergodic,
this implies that $D$ is of full measure in $S^1$.

Suppose for contradiction that the support $S$ of the measure
$\mu$ is not the whole $S^1\times\R$. Then for any $x\in D$
and $y\in\R$, the point $(x,y)$ cannot belong to $S$, for
if it did, then we would have  $S=S^1\times\R$ since the
orbit of $(x,y)$ is dense in $S^1\times\R$.
Therefore for the projection $q:S^1\times\R\to S^1$, we have
$q_*\mu(D)=0$. But $q_*\mu$ is just the Lebesgue measure $dx$,
contradicting the fact that $D$ is of full measure.
This completes the proof of Theorem \ref{T} (2).

\section{Appendix A: Properties of NICIS homeomorphisms}

Let $f$ be an area preserving homeomorphism of the annulus $\A$, 
isotopic to the identity.
Fix once and for all a lift $\tilde f$ of $f$ to the
universal covering space $\tilde\A$ of $\A$.\
As before the function $p\circ\tilde f-p$, defined on $\tilde \A$
and invariant by the deck transformation group, is considered
to be a function defined on $\A$. 

Denote by $\PP$
the space of the probability measures on $\A$ invariant by $f$.
Define
 $$
\mathcal R(\tilde f)=\{\langle p\circ\tilde f-p,\mu\rangle\mid
\mu\in\PP\}
$$ and call it the
{\em rotation set of} $\tilde f$.
The rotation set is a compact interval (maybe one point) of $\R$.

\begin{proposition} \label{p1}
If $f$ is a NICIS homeomorphism, then the rotation set
$\mathcal R(\tilde f)$ consists of a single irrational number $\alpha$.
\end{proposition}

{\sc Proof.}
It is shown in \cite{F1} that if $\mathcal R(\tilde f)$ is not a
singleton, 
then
there is a periodic orbit of $f$ in the interior $\IA$ of $\A$, contradicting
the condition of a NICIS homeomorphism. Assume that the rotation set
$\mathcal R(\tilde f)$ consists of a rational number $\alpha$. In the ergodic
decomposition of the Lebesgue measure, there must exist an
ergodic component $\nu$ whose support has nonempty intersection with
$\IA$. Since 
$$
\langle p\circ\tilde f-p,\nu\rangle=\alpha,
$$
 again there must
be a periodic orbit in $\IA$, by Proposition 2.3
of \cite{FH2}. 
\qed

\medskip
Next we shall show some topological properties of the
orbits of NICIS homeomorphisms.

\begin{proposition} \label{p2}
A NICIS homeomorphism $f$ admit dense orbits.
\end{proposition}

{\sc Proof.}
We shall show that any $f$-invariant open subset $U\subset\IA$ is dense
in $\IA$. By the Baire property of $\IA$, this will imply the existence of
dense orbits. In fact, choose a countable open basis $U_j$ of $\IA$.
The set $\mathcal O_f(U_j)$ 
of points whose orbits intersect $U_j$ will then be an open dense
subset of $\IA$, and a point in the nonempty intersection 
$\cap_j\mathcal O_f(U_j)$ will have a
dense orbit.

Suppose $U$ is an $f$-invariant open subset of $\IA$. Assume $U$ admits an inessential connected component $V$,
{\em i.\ e.} a component for which the inclusion $V\to\IA$ induces
the trivial map on the fundamental group. Then $V$ is left invariant
by some power of $f$, say $f^m$, since $f$ is area preserving.
Let $\hat V$ be the union of all the closed discs in $\IA$ whose boundary
is contained in $V$. Then $f^m$ leaves the open disc
$\hat V$ invariant, and admits
 a fixed point in $\hat V$ by the Brouwer plane fixed point theorem,
since all the points of $\hat V$ are nonwandering w.\ r.\ t.\
$f^m$. This contradiction
shows that any component of $U$ must be essential. 

Assume for contradiction that $U$ is not dense in $\IA$. By the same
argument applied to the interior of the complement of $U$, we obtain
two disjoint essential $f$-invariant open subset of $\IA$.
This implies the existence of interior compact invariant set.
A contradiction.
\qed

\medskip
Next proposition due to \cite{LY} and \cite{F2} 
shows that NICIS homeomorphisms admit nondense orbits as well.
Let us denote $\partial_{\pm}\A=S^1\times\{\pm 1\}$.

\begin{proposition} \label{p3}
Given any neighbourhood $U$ of $\partial_+\A$ (or of  $\partial_-\A$),
a NICIS homeomorphism $f$ admits an orbit contained in $U\cap\IA$.
\end{proposition}

\section{Appendix B: Proof of Theorem 3.1}

Throughout this section we assume that the skew product
$F_{\alpha,\varphi}$
satisfies the conditions (3.1) $\sim$ (3.4).

Computation shows that for any integer $n>0$,
$$
F_{\alpha,\varphi}^n(x,y)=(R^n_\alpha(x),y+\varphi_n(x)),
\ \ {\rm where} \ \
\varphi_n(x)=\sum_{i=0}^{n-1}\varphi(R^i_\alpha(x)).
$$

A positive integer $q$ is called a {\em closest return time} for
$\alpha$ if
$\Vert j\alpha \Vert > \Vert q\alpha \Vert$
for any integer $j$ such that $0<j<q$, where $\Vert\cdot\Vert$ denotes
the distance to 0 of the projected image in $\R/\Z$.
The following theorem can be found
in Sect.\ 1.2 of \cite{GLL}.

\begin{theorem} {\rm (Improved Denjoy-Koksma theorem)} \label{t3}
If $\varphi$ satisfies {\rm (3.2)} and {\rm (3.4)},
and $q_n$ is a sequence of closest return times for the
irrational number $\alpha$,
then we have
$$
\sup\abs{\varphi_{q_n}}
\to 0.
$$
\end{theorem}

It follows at once that for any $x\in S^1$, the closed subsemigroup
$\LL_x$ is nonempty. 

We also need the following lemma,
the proof of which is left to the reader.

\begin{lemma} \label{l1}
For integers $0\leq n_1< n_2< n_3$ assume that $n_3$
is the smallest positive integer $j$ such that $j \alpha$
lies in the smaller interval bounded by $n_1\alpha$ and
$n_2\alpha$.
Then $n_3-n_2$ is a closest return time for $\alpha$.
\qed
\end{lemma}

The following proposition plays a key role in the proof of Theorem 3.1.

\begin{proposition} \label{p6}
For any $x\in S^1$ and for any $0<a<b$, we have
$$\LL_x\cap([-b,-a]\cup[a,b])\neq\emptyset.$$
\end{proposition}

\smallskip

{\sc Proof.}
It is no loss of generality to assume that the point $x$ in Proposition
\ref{p6} is 0.
For contradiction assume that the limit set 
$\lambda(0,0)=\alpha(0,0)\cup\omega(0,0)$ 
is disjoint from
the set $0\times([-b,-a]\cup[a,b])$.
Thus one can choose a small positive
number $\varepsilon$ such that the union of two rectangles
$$
Y'=\{(x,y)\mid \vert x\vert\leq\varepsilon, \ a\leq\vert y\vert\leq b\}
$$
is disjoint from the orbit of $(0,0)$.

\smallskip
{\sc Claim.} {\em For even smaller $\delta>0$,
the orbit of $(0,0)$ is disjoint from
$\{\abs{x}\leq\delta,\abs{y}\geq a\}$.}

\smallskip

The claim clearly implies that the orbit of $(0,0)$ is
bounded, which contradicts the nonintegrability condition (3.3)
by Theorem 14.11 of \cite{GH}.

Let us start the proof of the claim. We are going to show
this only for the positive orbit, the other case being similar.

Choose positive integers $n_-$ and $n_+$ such that

\smallskip
\noindent
(5.1) $ -\varepsilon < n_- \alpha < 0 < n_+ \alpha < \varepsilon,$

\smallskip
\noindent
(5.2) $n_\pm$ {\em is a closest return time for $\alpha$, and}

\smallskip
\noindent
(5.3) {\em if $q$ is a closest return time for $\alpha$ such that
$q\geq\min\{n_-,n_+\}$, then}
 $$\sup\abs{\varphi_q}<b-a.$$

\smallskip
Put
$$
Y=\{(x,y)\mid n_-\alpha<x<n_+\alpha, \ a\leq\vert y\vert\leq b\}.
$$

Let us show that whenever a positive iterate of $(0,0)$ by $F_{\alpha,\varphi}$
lies in the strip $(n_-\alpha,n_+\alpha)\times \R$, it is contained
in $(n_-\alpha,n_+\alpha)\times(-a,a)$, that is, 
$n_-\alpha<n\alpha<n_+\alpha$ implies that $\vert \varphi_n(0)\vert<a$.
Assume to fix the idea that $n \alpha \in (0, n_+ \alpha)$.
Define a sequence $n_1, n_2,\cdots, n_k=n$ as follows.
First of all let $n_1=n_+$, and let $n_2$ be the
smallest positive integer $j$ such that $j \alpha$ lies in the interval
$I_1=(0, n_+ \alpha)$. If $n_2=n$, the definition is over. If not,
$n_2 \alpha$ divides the interval $I_1$ into two subintervals,
one of which, say $I_2$, contains $n \alpha$.
Let $n_3$ be the smallest positive integer $j$ such that $j \alpha$
lies in the interval $I_2$. Proceeding in this way, we end up with
some $n_k$ eventually matching the given $n$.

Now by Lemma \ref{l1}, all the integers 
$n_1, n_2-n_1, \cdots,n_k-n_{k-1}$ are closest return times for $\alpha$.
Moreover since $\Vert (n_j-n_{j-1})\alpha\Vert<\Vert n_+\alpha\Vert$,
we have $n_j-n_{j-1}\geq n_+$ ($2\leq j\leq k$) and (5.3) implies that
$\sup\abs{\varphi_{n_j-n_{j-1}}}\leq b-a$ as well as
$\sup\abs{\varphi_{n_1}}\leq b-a$.

To begin with, we have $\abs{\varphi_{n_1}(0)} \leq b-a<b$.
But since the point $F_{a,\varphi}^{n_1}(0,0)=(R_\alpha^{n_1}(0),
\varphi_{n_1}(0))$ does not lie in $Y$, we have in fact
$\vert \varphi_{n_1}(0)\vert < a$. 
Furthermore we obtain
$$
\abs{\varphi_{n_2}(0)}
=\abs{\varphi_{n_2-n_1}(R^{n_1}_\alpha(0))+\varphi_{n_1}(0)}
\leq
\abs{\varphi_{n_2-n_1}(R^{n_1}_\alpha(0))}+\abs{\varphi_{n_1}(0)}
\le
(b-a)+a.
$$
 Again since 
$F_{\alpha,\varphi}^{n_2}
(0,0)$ does not lie in $Y$, we actually have 
$\vert\varphi_{n_2}(0)\vert  <a$. 
Proceeding this way one gets $\vert
\varphi_{n}(0)\vert < a$.

The proof of the claim and of 
Proposition \ref{p6} is now complete.
\qed

\medskip
Let us finish the proof of Theorem \ref{t4}.
Proposition \ref{p6} shows that for any $x\in S^1$, the
set $\LL_x$ is unbounded. If the closed subsemigroup
$\LL_x$ contains both
positive and negative numbers, then $\LL_x=\R$ or $\LL_x=\Z\cdot c$
for some $c>0$. Proposition \ref{p6} shows that $\LL_x=\R$ in this 
case. Also the same proposition
shows that if $\LL_x$ consists of nonnegative numbers
(resp.\ nonpositive numbers), then $\LL_x=[0,\infty)$
(resp.\ $\LL_x=(-\infty,0]$).
This shows that $S^1$ is the union of
$D$, $P$ and $N$.

What is left is to show that
each of the sets $D$, $P$ and $N$ is nonempty. 
The existence of dense orbits of $F_{\alpha,\varphi}$
is shown in \cite{Hed} under the assumption (3.1) $\sim$ (3.3), which implies
that $D$ is nonempty.
 In fact the essential point is to show that
any point of $S^1\times\R$ is nonwandering, the rest being
exactly the same as the proof of Proposition \ref{p2}.
Also in \cite{B} the existence of the $F_{\alpha,\varphi}$-orbits
which are bounded below (or above) is shown under the same assumption.
This, together with Proposition \ref{p6},
shows that $P$ and $N$ are nonempty.


\begin{thebibliography}{99}


\bibitem[B]{B}
A. S. Besikovich, {\em A problem on topological transformations of the
plane II,} Proc.\ Cambridge Phi.\ Soc.\ {\bf 47}(1951), 38-45.


\bibitem[F1]{F1}
J. Franks, {\em Generalization of the Poincar\'e-Birkhoff theorem,}
Ann.\ Math.\ {\bf 128}(1988), 139-151.

\bibitem[F2]{F2}
J. Franks, {\em The Conley index and non-existence of minimal 
homeomorphisms,} Illinois J. Math.\ {\bf 43}(1999), 457-464.

\bibitem[FH2]{FH2}
J. Franks and M. Handel, {\em Distortion elements in group actions
on surfaces,} Duke Math.\ J. {\bf 131}(2006), 441-468.

\bibitem[FS]{FS} B. Fayad and M. Saprykina, {\em Weak mixing disc and
	annulus
diffeomorphisms with arbitrary Lioubille rotation number on the
	boundary},
Ann.\ Sci.\ \'Ecole Norm.\ Sup.\ (4) {\bf 38} (2005) no.\ 3, 339-364.

\bibitem[GH]{GH} W. H. Gottschalk and G. A. Hedlund,
{\em Topological Dynamics,} Amer.\ Math.\ Soc.\ Colloquium Publication
Vol.\ {\bf 36}(1955).

\bibitem[GLL]{GLL} P. Gabriel, M. Lemanczyk and P. Liardet,
{\em Ensemble d'invarinats pour les produits crois\'es de Anzai,}
Mem.\ Soc.\ Math.\ France {\bf 47}(1991), 1-102.




\bibitem[Hed]{Hed}
G. A. Hedlund, {\em A class of transformations of the plane,}
Proc.\ Cambridge Phil.\ Soc.\ {\bf 51}(1955), 554-564.

\bibitem[He1]{He1}
M. R. Herman, {Some open problems in dynamical systems,} Proc.\ ICM,
Berlin 1998, Vol.\ 2, p.797-808.


\bibitem[HR]{HR}
G. H. Hardy and W. W. Rogosinski, {\em Fourier series},
Cambridge Tracts in Mathematics and Mathematical Physics,
No.\ 38 (1950).



\bibitem[K]{K}
A. B. Krygin, {\em $\omega$-limit set of smooth cylindrical cascades,}
Math.\ Notes {\bf 23}(1978), 479-485.


\bibitem[LY]{LY}
P. Le Calvez and J.-C. Yoccoz, {\em Un th\'eor\`em d'indice pour
les hom\'eomorphismes du plan au vousinage d'un point fixe,}
Ann.\ Math.\ {\bf 146}(1997), 241-293.

\bibitem[M]{M} J. Moser, {\em On the volume elements on a manifold,}
Trans.\ Amer.\ Math.\ Soc.\ {\bf 120}(1965), 286-294.


\bibitem[OU]{OU}
J. Oxtoby and S. Ulam, {\em Measure preserving homeomorphisms and
metrical transitivity,}
Ann.\ Math.\ {\bf 42}(1941), 874-920.

\end{thebibliography}
\end{document}